\documentclass[]{interact}

\usepackage{epstopdf}% To incorporate .eps illustrations using PDFLaTeX, etc.
\usepackage[caption=false]{subfig}% Support for small, `sub' figures and tables
\usepackage[natbibapa,nodoi]{apacite}% Citation support using apacite.sty. Commands using natbib.sty MUST be deactivated first!
\usepackage[colorlinks=true,citecolor=black,linkcolor=black]{hyperref}

\theoremstyle{plain}% Theorem-like structures provided by amsthm.sty

\theoremstyle{definition}

\theoremstyle{remark}

\newcommand{\citeshortp}[1]{\shortcites{#1}\citep{#1}}
\newcommand{\citeshortt}[1]{\shortcites{#1}\citet{#1}}

\begin{document}

\noindent
The Version of Record of this manuscript has been published and is available in Structure and Infrastructure Engineering, 16 January 2024, \url{http://www.tandfonline.com/10.1080/15732479.2023.2297891}.
\newpage

% \articletype{ARTICLE TEMPLATE}% Specify the article type or omit as appropriate

\title{Socio-technical systems integration and design: a multi-objective optimisation method based on integrative preference maximisation}

\author{
\name{Harold van Heukelum\textsuperscript{a}*\thanks{*Corresponding author. Email: harold.van.heukelum@boskalis.com.
% \\
% The AOM of this paper has been published on the pre-print server of ArXiv (Van Heukelum, Binnekamp, et al., \protect\citeyear{vanheukelum2023human}).
}, Ruud Binnekamp\textsuperscript{b} and Rogier Wolfert\textsuperscript{c}}
\affil{\textsuperscript{a}Corporate R\&D, Royal Boskalis, Rosmolenweg 20, 3356 LK Papendrecht, The Netherlands; \textsuperscript{b}Department Materials, Mechanics, Management \& Design, Faculty of Civil Engineering and Geosciences, Delft University of Technology, Stevinweg 1, 2628 CN Delft, The Netherlands; \textsuperscript{c}Department of Engineering Structures, Faculty of Civil Engineering and Geosciences, Delft University of Technology, Stevinweg 1, 2628 CN Delft, The Netherlands}
}

\maketitle

\begin{abstract}
Current systems design optimisation methodologies are one-sided, as these ignore the socio-technical integration between stakeholder preferences (‘what a human wants’) and the capability of technical assets (‘what a system can deliver’). Moreover, classical multi-objective optimisation methods contain fundamental mathematical flaws. Also, the often-used classical Pareto front does not provide a single best-fit design configuration, but rather a set of design alternatives. This leaves designers without a unique solution to their problems. Finally, current multi-objective optimisation processes are not well aligned with design practices, because they do not sufficiently involve decision makers and do not translate their interests into a single common preference domain to find an overall group optimum. This paper introduces a new Open Design Systems (Odesys) methodology and a new Integrative Maximisation of Aggregated Preferences (IMAP) method, implemented in the Preferendus tool. Its added value and use are exemplified in two infrastructure design applications, which show how to achieve the pure best-fit for common-purpose design results.
\end{abstract}


\section{Introduction}
Why, so often, do people build what nobody wants? Why, so often, do engineers optimise their solution based only on physical capabilities and fail to consider the stakeholders’ desirabilities? Why, so often, do infrastructure managers keep the design/decision-making process non-transparent and non-participatory? The answer to these questions is that engineering design and decision-making are often solved from a one-sided point of view, without considering the fact that the problem is complex and multifaceted. Therefore, a participatory process that does justice to both the ‘hard’ technical and ‘soft’ social aspects of infrastructure systems development is needed. It is thus crucial to truly connect and bridge the gap between human preference interests and the engineering assets performances using transparent models for complex systems design and integration. The goal of such an Open Design Systems (Odesys) approach is to promote the use of the civil infrastructures that surround us every day through a multi-system level socio-technical approach, supported by sound mathematical open-glass box models as a means of observation and perception in collaborative decision-making.

Above all, zooming in on the design challenge of our contemporary (civil) infrastructures, it can be noted that this challenge is becoming increasingly complex due to environmental demands, new transport modes, and other transitions. This rapidly changing infrastructure context requires an optimal life-cycle value design within the framework of infrastructure asset management (see e.g. \citet{hasting2015,balzer2016}; \citeshortt{uddin2013public} or the NEN 15288 on system life-cycle processes (\citeyear{nen15288})). Multi-objective optimisation is key to supporting informed decisions in infrastructure asset management (for an extensive literature review and overview of optimisation methods, see \citet{chen2019optimization}). Increasing stakeholder involvement, combined with the multidisciplinary nature of infrastructure design challenges, further necessitates a more effective and efficient participatory and supportive decision-making process (among others, see also \citeshortt{omar2009infrastructure,wang2020approaches,chen2021multi}). In this context of asset management decision-making, the focus of this paper is therefore on socio-technical design optimisation, where both the various stakeholder preference interests (or societal values) and the technical system life-cycle capabilities are unified in a best-fit for common-purpose design configuration. To this end, a new so-called Odesys methodology is introduced, with a new preference-based multi-objective design optimisation method. This is required because the current design optimisation methods have intrinsic problems and/or shortcomings that make them unsuitable to provide the required unique and best-fit design solutions. These five fundamental problems and shortcomings together constitute Odesys’ \emph{development gap}.

The first problem with the current multi-objective design optimisation methodologies is the disconnect between the domain of human preferences (subject desirability) and the domain of the physical performance behaviour of the engineering asset (object capability). Moreover, when applied in the classical systems engineering context, design optimisation is usually limited to a single objective design approach and/or to an \textit{a posteriori} evaluation of design alternatives \citep{Dym2004,Blanchard2011,cross2021engineering}. However, in \textit{a posteriori} evaluation, there is no guarantee that the optimal design point has been found and a choice has to be made between sub-optimal compromise solutions (even when optimisation and \textit{a posteriori} evaluation are combined, see \citet{mueller2015combining}). Especially in complex engineering projects, the number of possible design alternatives is too large to evaluate them all and the optimal solution may thus be ignored.

Secondly, most multi-objective optimisation methodologies introduce fundamental mathematical operation and aggregation flaws because they: (1) use undefined measurement scales and apply mathematical operations where these are not defined (e.g. for variables that have neither an absolute zero nor one, such as time/potential energy/preference, the mathematical operations of addition and multiplication are not defined in the corresponding mathematical model which is the one-dimensional affine space); (2) produce an infinite number of non-equivalent ‘optimal’ outcomes (e.g. the definition of the aggregation algorithm does not prerequisite having only normalised numbers); (3) outcomes do not take into account the relative scoring impact of other design alternatives (e.g. in reality, the score of one alternative depends on the performance of all the other alternatives; the score is obtained by finding the best balance between the normalised and weighted scores for all sub-criteria given the set of alternatives). As a result, the outcomes of decision-making in engineering design may lead to sub-optimal design configurations. The foundations of this second shortcoming are found from the principles of Barziali’s Preference Function Modelling (PFM) and its associated preference measurement theory \citep{barzilai2005measurement,barzilaiEngineeringDesign,barzilai2022pure}.

A third problem with many of the classical multi-objective design optimisation methods is that they do not have a consistent way of translating the different objective functions into a common domain to find a best-fitting aggregated optimum. To get around this problem, these multi-objective design methods often use monetisation. In other words, all objective functions are expressed in terms of money. However, according to classical decision/utility theory, decisions are not based on money, but on value or preference (where minimising expenditure or maximising profit can be one of the objectives). Here, preference is an expression of the degree of ’satisfaction’, and it describes the utility or value that something provides. Although some researchers have incorporated preference modelling into their multi-objective optimisation frameworks (see, for example, \citeshortt{lee2011preference} or \citet{messac1996physical}), none of them use strong (preference) measurement scales or individually weighted preference functions (i.e. continuous functions linking an individually weighted preference to a specific objective). In addition, these approaches do not lead to a single optimal design point and also contain the aggregation modelling errors mentioned above.

A fourth shortcoming of classical multi-objective design optimisation methods is that many of them consider the so-called Pareto front as a valid outcome \citep{marler2004survey}. Apart from the fact that the Pareto front is often obtained in a mathematically incorrect way (see the aforementioned second point), it also generates an infinite set of possible, and supposedly equally desirable, design points (see, for example, \citeshortt{farran2015fitness,furuta2006optimal,saad2018concurrent}). However, this is inconsistent with the fundamental basis of an engineering design process, where each design point is (subjectively) interpreted by people in terms of preference (i.e. a statement of their individual interest) and where a search is performed to find a single optimal design solution. These Pareto shortcomings are also noted by e.g. \citeshortt{kim2022probabilistic,lee2011preference,bai2015,golany2006efficiency,bakhshipour2021toward}, amongst others. However, their proposed (hybrid) solutions still rely on the Pareto front (with its mathematical flaws) and some form of \textit{a posteriori} evaluation. Their modelling approaches therefore fail to provide a pure integrative design approach and are not able to obtain \textit{a priori} a single best configuration.

A fifth shortcoming is that current multi-objective optimisation processes are rather disconnected from systems design practices, as they lack deep involvement of decision-making stakeholders \citep{guo2022multi}. In addition, the dynamic nature and the socio-technical interaction between stakeholder preferences (‘what a human wants’) and the performance of technical assets (‘what a system can’) are often not considered in service life design.

To overcome the aforementioned shortcomings and problems, and to enable pure human preference and asset performance systems design integration, the Odesys design methodology is introduced in this paper. Odesys builds further on the multi-stakeholder design optimisation methodology proposed by \citeshortt{zhilyaev2022best}, who showed that the unambiguous solution to a multi-objective engineering design/decision problem is to translate each of the objective functions, as a function of the design variables, into an overarching preference domain. This can be done using stakeholder preference functions: i.e. the relationship between an individual preference and a specific objective, which then allow for the maximisation of the aggregated group preference, leveraging Barzilai’s PFM theory (see \citet{binnekamp2010preference} where this concept originated in its initial form, and \citeshortt{arkesteijn2017improving} for its early social validation). However, all these aforementioned developments in the field of preference-based design, which so far only were applied in the context of real estate planning, still have three methodological deficiencies, and lack the following:
\begin{enumerate}
    \item a generalised mathematical framework for multi-objective socio-technical design optimisation: i.e. a threefold modelling framework of integrative performance, objective and preference functions;
    \item a connection between common socio-eco interests and the related subject preferences, and the physical/mechanical object behaviour: i.e.  a pure integration of technical design performance, social objective and preference functions;
    \item a PFM-based solver: i.e. a search algorithm to find the optimal solution with the maximum aggregated preference.
\end{enumerate}

\noindent
As a conclusion to the five fundamental shortcomings of current multi-objective design optimisation methods and the three deficiencies of the preference-based design approach mentioned above, the Odesys \emph{development statement} reads as follows:\\

\emph{“There is a need for an open design/decision methodology enabling socio-technical systems integration on all relevant levels using a human-centred preference-based design performance approach supported by pure mathematical optimisation modelling.”}\\

\begin{figure}
    \centering
    \includegraphics[width=\textwidth]{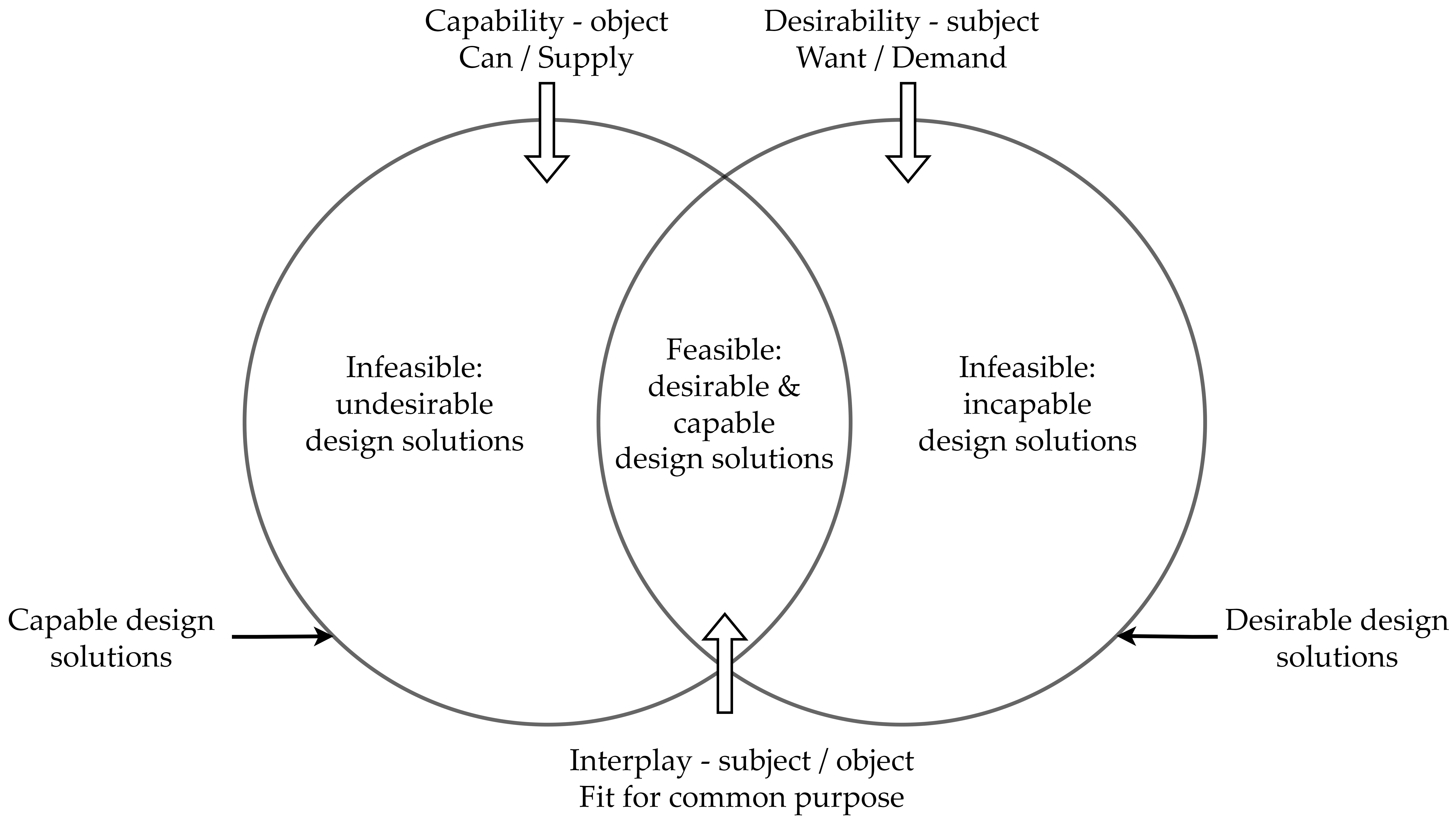}
    \caption{Socio-technical interplay between (un)desirability and (in)capability.}
    \label{fig:venn_diagram}
\end{figure}

This is the basis for the development of the Odesys methodology methods and tool, which will allow for the full integration between subject (un)desirability: ‘what a stakeholder wants/does not want’, expressed via preference functions, and object (in)capability: ‘what a system can/cannot’, expressed via design performance functions. This integration is schematically depicted in \autoref{fig:venn_diagram}. It is being achieved by constructing preference functions that are a direct function of both the stakeholder objective and the engineering asset design performance functions, which depend on the design and physical variables and their constraints. In other words, this unified set of preference functions, which at the lowest level is a function of the engineering design variables and the physical constraints, is a translation (a mapping) of the socio-technical system under consideration. Next, an automated algorithm is needed that searches for a feasible and optimal design synthesis solution where the aggregated group preference score is maximal. In reality, this search is an open-ended approach. This means that an iterative process of technical-, social-, and purpose-cycles will have to take place. This implies that a best-fit for common-purpose design configuration can only be achieved through an iterative socio-technical process given the final ‘idealised’ desires, objectives, interests, and requirements of the stakeholders.

This makes Odesys a pure socio-technical systems integration methodology where human preference-based design and engineering physics/mechanics converge, offering a wide range of potential applications within the context of (infra)structure systems engineering design. As part of this Odesys methodology, a new Integrative Maximised Aggregated Preference (IMAP) optimisation method for maximising aggregated preferences is introduced. This IMAP method forms the basis of a new software tool called the Preferendus and combines the state-of-the-art PFM principles with an inter-generational Genetic Algorithm (GA) solver developed specifically for this purpose. It should be noted that the Preferendus in its primary form was published by \citeshortt{zhilyaev2022best}.

This paper continues  by giving a general mathematical statement of the Odesys methodology. Next, a flow chart (or concept diagram) of the Preferendus software tool is described, in which the Odesys methodology is implemented. Finally, the use and added value of the Odesys methodology, the IMAP optimisation method, and the Preferendus tool are demonstrated for two infrastructure life-cycle design applications carried out in a real industrial context (a marine contractor and a railway infrastructure service provider). The results are compared with single-objective design outcomes as well as with design outcomes resulting from the classical min-max goal attainment method.

\section{Mathematical formulation of the Open Design Systems methodology}
As described in the introduction, there is currently no optimisation framework that allows for pure integration of the human preference domain (subject desirability) and the engineering asset physical performance behaviour domain (object capability). This disconnection will limit optimisation to sub-optimal results, as the interaction between these two levels is not considered. To overcome this, this paper introduces the following mathematical statement, which integrates subject desirability and object capability, and is at the core of the Odesys methodology:

\begin{equation}
    \label{eq:general_MS}
    \begin{gathered}
    \mathop{Maximise}_{\mathbf{x}} \textrm{ } U = T \left[  
    P_{k,i} 
    \left(O_i
    \left(F_1(\mathbf{x}, \mathbf{y}),F_2(\mathbf{x}, \mathbf{y}),...,F_J(\mathbf{x}, \mathbf{y}) \right)
    \right),w'_{k,i} \right] \textrm{ for } \\ 
    k=1,2,...,K \\ 
    i=1,2,...,I
    \end{gathered}
\end{equation}

Subject to: 
\begin{equation}
    \label{eq:ineq_cons}
    \begin{gathered}
    g_{p}(O_{i}(F_{1,2,...,J}(\mathbf{x}, \mathbf{y})),F_{1,2,...,J}(\mathbf{x}, \mathbf{y})) \le 0
    \textrm{ for } p=1,2,...,P
    \end{gathered}
\end{equation}

\begin{equation}
    \label{eq:eq_cons}
    \begin{gathered}
    h_{q}(O_{i}(F_{1,2,...,J}(\mathbf{x}, \mathbf{y})),F_{1,2,...,J}(\mathbf{x}, \mathbf{y})) = 0
    \textrm{ for } q=1,2,...,Q
    \end{gathered}
\end{equation}

\noindent
With:

\begin{itemize}
    \item $T$: The aggregated preference score determined using the PFM theory principles (see \citet{barzilai2022pure}).
    
    \item $P_{k,i}(O_i(F_{1,2,...,J}(\mathbf{x}, \mathbf{y})))$: Preference functions that describe the preference stakeholder $k$ has towards objective functions, which are functions of different design performance functions and dependent on design and physical variables.
    
    \item $O_i(F_{1,2,...,J}(\mathbf{x}, \mathbf{y}))$: Objective functions that describes the objective $i$, functions of different design performance functions and dependent on design and physical variables.
    
    \item $F_{1,2,...,J}(\mathbf{x}, \mathbf{y})$: Design performance functions that describe the object, depending on one or multiple design variables x (i.e. controllable endogenous variables) and one or multiple physical variables y (i.e. uncontrollable exogenous variables).
    
    \item $\mathbf{x}$: A vector containing the (controllable) design variables $x_1,x_2,...,x_N$. These variables are bounded such that $lb_n \le x_n \le ub_n$, where $lb_n$ is the lower bound, $ub_n$ is the upper bound, and $n=1,2,...,N$.
    
    \item $\mathbf{y}$: A vector containing the (uncontrollable) physical variables $y_1,y_2,...,y_M$.
    
    \item $w'_{k,i}$: Weights for each of the preference functions. These weights can be broken down into weights for the stakeholders and weights for the objectives:
    \begin{itemize}
        \item $w_{k}$: weights for stakeholders $k=1,2,...,K$. These weights represent the relative importance of stakeholders.
        \item $w_{k,i}$: these weights represent the weight stakeholder $k$ gives to objective $i$.
    \end{itemize}
    
    \noindent
    The final weights $w'_{k,i}$ can be constructed via $w'_{k,i}=w_{k}\cdot w_{k,i}$, given that $\sum w'_{k,i}=\sum w_{k,i}=\sum w_k=1$

    \item $g_{p}(O_{i}(F_{1,2,...,J}(\mathbf{x}, \mathbf{y})),F_{1,2,...,J}(\mathbf{x}, \mathbf{y}))$: Inequality constraint functions, which can be either objective function and/or design performance function constraints.
    \item $h_{q}(O_{i}(F_{1,2,...,J}(\mathbf{x}, \mathbf{y})),F_{1,2,...,J}(\mathbf{x}, \mathbf{y}))$: Equality constraint functions, which can be either objective function and/or design performance function constraints.
\end{itemize}

\noindent
To further elaborate on this formulation, several important remarks are made which are discussed below.

\subsubsection*{Remark 1: preference aggregation}
Here, the aggregated preference scores are determined based on the principles of PFM, expressed by the mathematical operator $T$. This operator is a solving algorithm that is based on finding/synthesising the aggregated preference score (i.e. the ‘best’ fit of all weighted (relative) scores for all the decision-making stakeholders’ objectives) that minimises the least-squares difference between this overall preference score and each of the normalised individual scores (on all criteria) by computing its closest counterpart \citep{barzilai2022pure,zhilyaev2022best}.

In this, preference is a statement of an individual stakeholder’s interest and a measure of satisfaction, which is a score that is expressed as a real number (scalar or bare quantity) on a defined scale, e.g. 0 to 100, where 0 corresponds to the ‘worst’ performing alternative and 100 to the ‘best’ performing alternative. For the applications shown in this paper, Tetra is used as this preference aggregation solver. For more information on the Tetra solver, see \citet{scientific_metrics}.

\subsubsection*{Remark 2: preference functions}
Preference functions describe the relationship between an individual stakeholder’s preference and a specific objective (where a stakeholder is defined as one of the participants in the design/decision-making process). The theory of preference functions (often also called utility functions) for \textit{a posteriori} multi-criteria decision evaluation is a branch of the social science in itself. However, the preference functions are needed as input to the design/decision system to enable \textit{a priori} multi-objective design optimisation.  Here, the elicitation of the preference functions and associated weights is handled pragmatically using ‘static’ expert judgement, whereas in practice this is inherently a dynamic and iterative process that helps stakeholders better understand the impact of their input on the optimisation outcome (see \citeshortt{arkesteijn2017improving} for the specifics of this elicitation as part of the design cycle).

Finally, note that an objective $O_i$ can be associated with multiple stakeholder preference functions $P_{k,i}$ (as $k\geq i$). However, it is not required that a stakeholder expresses a preference for all objectives. This is modelled by giving a stakeholder’s objective a weight of zero, which means that some elements of the $w_{k,i}$ matrix can be zero.

\subsubsection*{Remark 3: bound preference scores}
Here, a preference score is bounded by $0\le P_{k,i}\le100$. A constraint can be added to the objective functions to prevent preference scores which lay outside these bounds.

\subsubsection*{Remark 4: design variables in objective functions}
A design variable $x$ can be directly linked to an objective function $O$. In this case, the design performance function $F$ is just equal to the design variable $x$. Moreover, these design performance functions $F$ can also only relate to an exogenous physical variable $y$.

\subsubsection*{Remark 5: rewrite equality constraints}
Equality constraints are quite common in the object behaviour domain. However, as the Preferendus uses a GA, equality constraints can complicate the convergence of the optimisation, as especially the simpler constraint handlers for GAs have problems with handling equality constraints \citeshortp{Kramer2017,homaifar1994constrained}. Therefore, when modelling a system of interest, the equality constraints can be rewritten as inequality constraints, as is often done in literature \citep{coello2002theoretical,Kramer2017}. This is often done in the form of Equation (\ref{eq:rewrite_eq_cons}).

For the proposed Odesys methodology, it is possible to rewrite most equality constraints directly into inequality constraints, as the methodology aims to reduce ‘waste’ in the result. For example, the length of a beam supporting a floor will usually have a fixed length: the length of the span. Since a length greater than the length of the span will result in more costs, material consumption, carbon emissions, etc., this equality constraint can safely be rewritten as an inequality constraint. This makes modelling easier, since the tolerance $\epsilon$ does not have to be set and tuned for each problem.

\begin{equation}
    \label{eq:rewrite_eq_cons}
    |h_{1,2,...,Q}(O_{1,2,...,I}(F_{1,2,...,J}(\mathbf{x},\mathbf{y})),F_{1,2,...,J}(\mathbf{x},\mathbf{y}))|-\epsilon \le 0
\end{equation}

\subsubsection*{Remark 6: soft and hard constraints}
Finally, a distinction can be made between soft and hard constraints. The former result from the sociological aspect of a design process and are negotiable. They can be adapted during the process based on discussions with other stakeholders or new insights. The latter are fixed and non-negotiable. They are given by, among others, laws of nature, material composition, environmental conditions, etc.

\subsection{Conceptual threefold framework}
The mathematical statement with the aforementioned remarks provides a general framework in which it is possible to connect the subject desirability level (preference functions) with the object capability level (design performance functions) via the integrative subject-object conciliation level (objective functions). Note that there will be three types of functional values and/or outcomes of interest: (1) degrees of capability – design performances (technical); (2) degrees of freedom – design variables (technical); (3) degrees of satisfaction – preferences and objectives (social).

\begin{figure}
\centering
\includegraphics[width=0.7\linewidth]{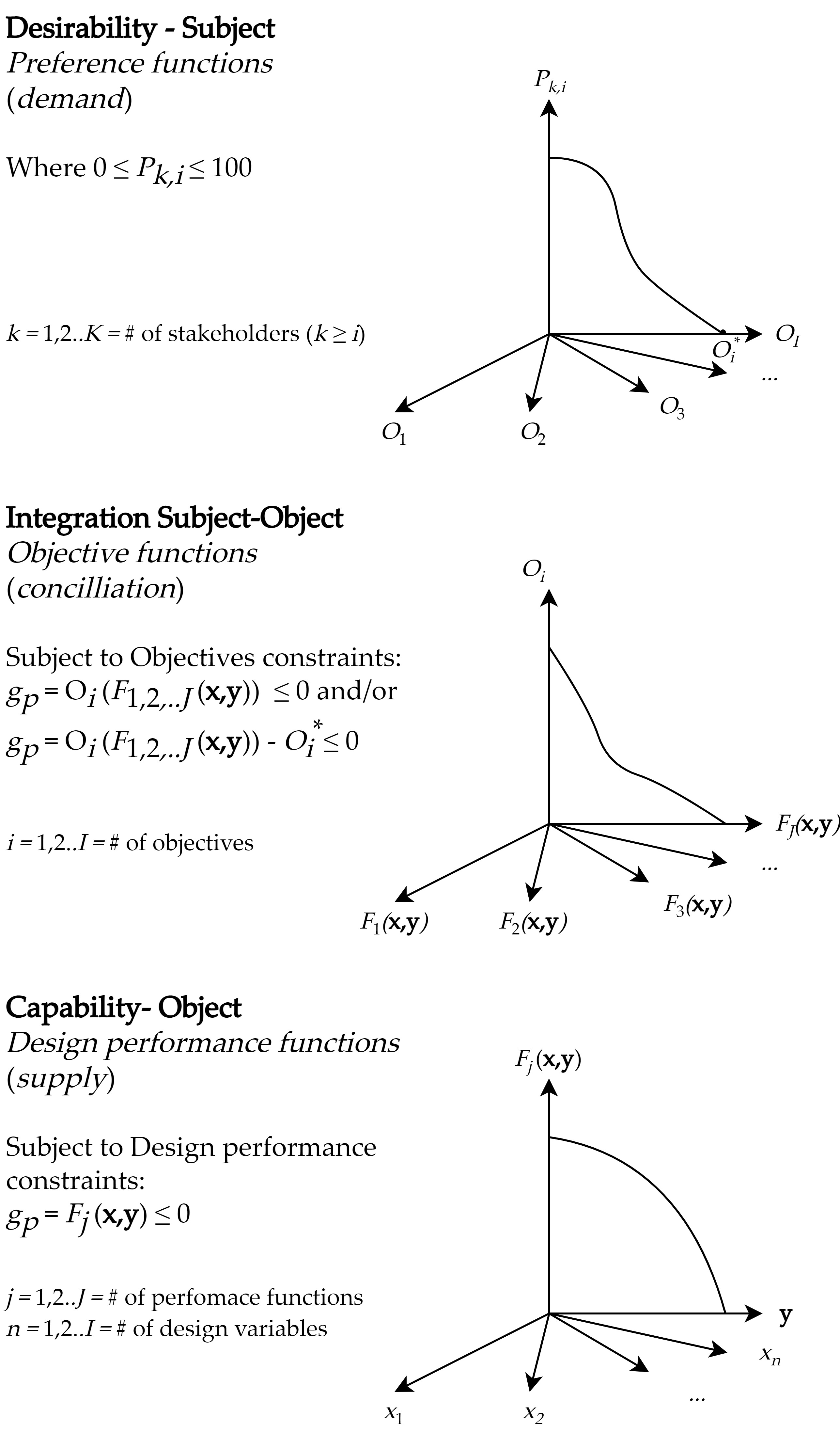}
\caption{Conceptual threefold framework of the Odesys mathematical statement, where subject desirability (preference functions) and the object capability (design performance functions) are integrated subject-object (objective functions). Note: the shapes of the curves are arbitrary.}
\label{fig:threefold}
\end{figure}

To better understand and further detail this specific social-technical systems integration, the different functions as part of the mathematical formulation are conceptualised in a threefold modelling framework, as shown in \autoref{fig:threefold}. Note that the different functions are linked (an ordering principle) and that maximisation is not yet part of this threefold.

\section{The Preferendus \& IMAP}
In this section, the so-called Preferendus tool and its optimisation method IMAP are described as part of the Odesys methodology. Here, the conceptual functioning will be introduced (as an extension and further development of the Preferendus as described by \citet{zhilyaev2022best}) and in the following section, design applications using the state-of-the-art Preferendus tool are demonstrated.

The Preferendus tool is based on the IMAP optimisation method presented in this paper. It combines proper preference aggregation with preference maximisation, as described by the mathematical formulation of the Odesys problem statement in the previous section.

\subsection{Preference aggregation (IMAP part 1: synthesis)}
Following the Odesys methodology, it is argued that the overarching goal of multi-objective design optimisation is to find the highest overall group preference score that represents the design synthesis. However, for these design syntheses to be possible, the individual preference scores first need to be aggregated.

Since preference scores are defined in an affine space, aggregation should also take place in this space. This means that, according to the basic principles of PFM theory, the correct way of aggregating preference scores is to find the aggregated preference score that provides the ‘best’ fit to all the weighted (relative) scores of the different preference functions ($P_{k,i}$). Here, the preference functions are the integration of objective functions and design performance functions. The final preference score aggregation is performed by the aforementioned PFM-based solving approach (see remark 1 of the previous section), as an integral part of the overall design optimisation algorithm.

\subsection{Preference maximisation (IMAP part 2: synthesis)}
To finally find the design configuration that reflects the maximum group preference aggregation, it is also necessary to use a maximisation algorithm. To do this, a GA is used that is specifically adapted to work with Tetra. This is necessary because it is not possible to directly compare one generation of the GA with another, as the aggregated preference scores contain only information about the alternatives of a single generation. To overcome this, a GA is developed that combines widely available elements and is extended with a so-called inter-generational solver. The details and the operation of this GA solver are given in \autoref{app:solver}.

\begin{figure}
\centering
\includegraphics[width=\linewidth]{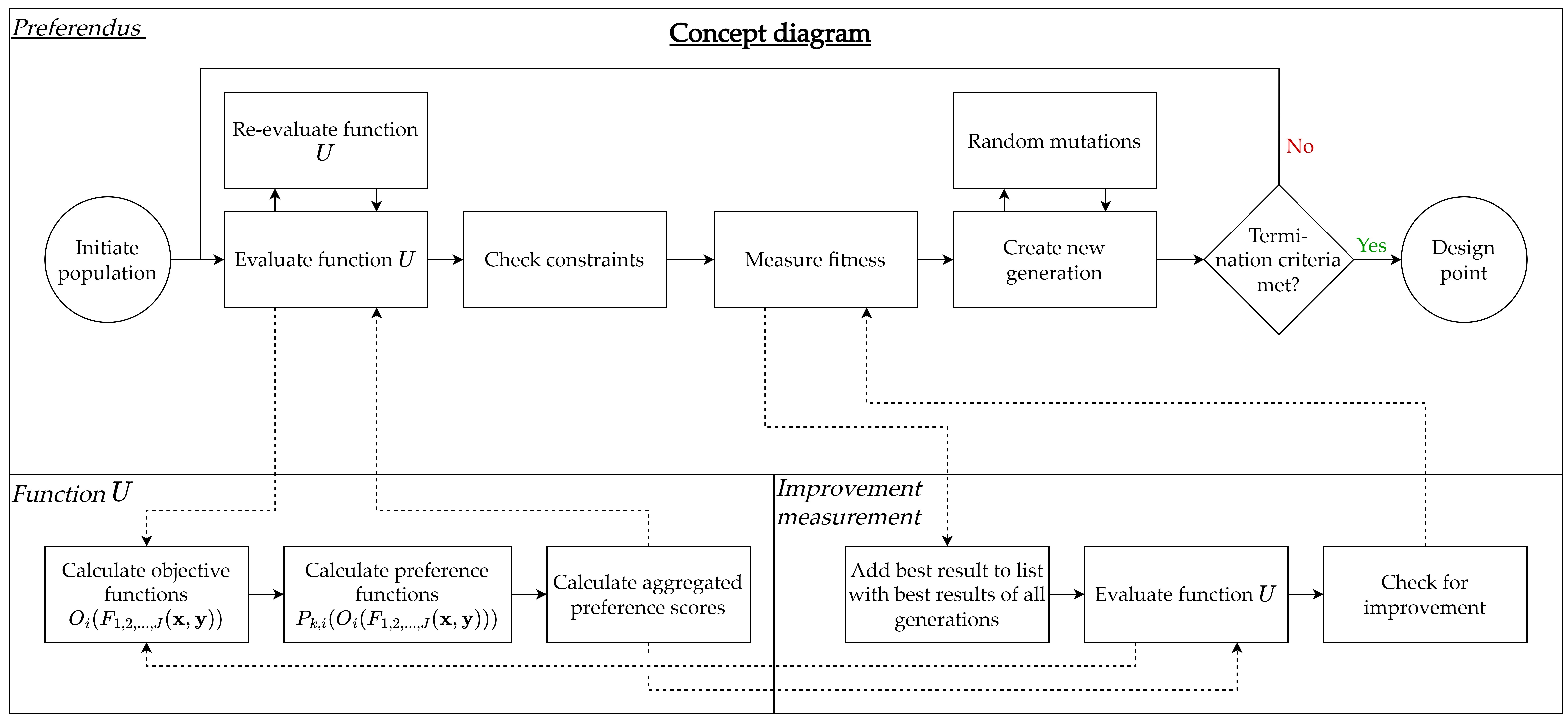}
\caption{The workflow of the Preferendus, presented as a concept diagram.}
\label{fig:concept_diagram_preferendus}
\end{figure}

The final result is an Odesys-based design optimisation tool, the Preferendus, which incorporates the IMAP method. The concept diagram of the Preferendus is shown in \autoref{fig:concept_diagram_preferendus}. This is an open-source tool available via GitHub (see data availability statement).

\subsection{Interlude min-max goal attainment (min-max: compromise)}
To compare and validate the results and added value of the IMAP multi-objective optimisation method, the following section first compares the results with those of the single-objective optimisation. In addition, a comparison is made with the classical min-max goal attainment multi-objective optimisation method \citep{marler2004survey}. This method does not generate group results based on overall aggregation, but rather optimises, i.e. equalises, each individual result so that it is as close as possible to a ‘utopian’ design point. In other words, the min-max method tries to minimise the maximum dissatisfaction for all individual scores (expressed by the distance to this utopia point). The result of this method, which does not conflict with the fundamental PFM principles, is a solution that gratifies each stakeholder equally.

In order to make a like-for-like comparison between IMAP and min-max, the mathematical formulation of the Odesys problem statement needs to be modified (i.e. Equation (\ref{eq:general_MS}) needs to be changed). First, this means that in this case the min-max method will try to minimise the distance to a score of 100 for all different preference scores $P_{k,i}$ (i.e. the best-scoring utopian point has been defined as 100). Then, the preference score $P_{k,i}$ with the greatest (weighted) dissatisfaction must be found and minimised, which mathematically can be read as Equation (\ref{eq:min_max}).

\begin{equation}
    \label{eq:min_max}
    \begin{gathered}
    \mathop{Minimise}_{\mathbf{x}} \textrm{ } U = \max_{k,i} \left[ 
    w'_{k,i} \times \left\{100 - P_{k,i}\left(O_i \left(F_1(\mathbf{x}, \mathbf{y}),F_2(\mathbf{x}, \mathbf{y}),...,F_J(\mathbf{x}, \mathbf{y}) \right) \right) \right\} \right] \textrm{ for } \\ 
    k=1,2,...,K \\ 
    i=1,2,...,I
    \end{gathered}
\end{equation}

\noindent
It should be noted that the min-max goal attainment method, as part of a larger group of multi-objective optimisation methods, does not violate the PFM principles. However, this method treats the scores of all design alternatives as absolute values, ignoring the dynamic interplay between them. In other words, this method focuses on making each stakeholder as ‘happy’ as possible, even though this may not be beneficial for the group as a whole. This is why this optimisation is called a compromise method, because it finds a design configuration based on a compromise between stakeholders rather than a synthesis.

\section{Real-life service life design applications}
The Odesys methodology, the associated IMAP optimisation method and the use of the Preferendus are demonstrated in two real-life infrastructure design applications: (1) a railway level-crossing life-cycle design and (2) a floating wind turbine installation design. The source code and results of these applications are available on GitHub, see the data availability statement.

Both of these design applications were conducted within a real-life infrastructure design context. The first in collaboration with ProRail, a Dutch railway infrastructure service provider, and the second within Boskalis, an internationally operating maritime contractor. Especially within Boskalis, several socio-technical cycles were carried out to validate the Odesys results and the added value of the Preferendus with various stakeholders involved. In addition, here the Preferendus has also been validated for a dredging application with promising validation results (but beyond the scope of this paper).

Although both design applications are simplified for illustrative purposes, they still provide insight into the added value and principles of the Odesys methodology, the IMAP method and the use of the Preferendus tool. For further substantial extensions of both design application cases as they are presented here, see Shang, Binnekamp, et al. (\citeyear{shangIALCCE}) and \citeshortt{vanheukelumIALCCE}.

For both design applications, the threefold diagram of design performance, objective, and preference functions is presented. Preference functions are ultimately a direct function of both stakeholder objective and engineering asset design performance functions, which in turn are related to the design variables and their constraints. These functions are derived from an idealised design configuration (i.e. a tangible design representation) and from the common preference interests of the stakeholders involved. The goal is then to find, within the feasibility space, the candidate solution with the highest aggregated group preference.

In real-life design practices, this quest is an open-ended approach. This means that an iterative process of technical-, social-, and purpose-cycles will have to take place, implying that a best-fit for common-purpose design configuration can only be achieved through an iterative socio-technical process given the final ‘idealised’ desires, objectives, interests, and requirements of the participating stakeholders. In this paper, for demonstration purposes, only one socio-technical cycle per design application is included. It should be noted that in the real-life design application 2 (‘floating wind installation’), as carried out within Boskalis, the stakeholders were asked to adjust their preference and/or objective functions (social context) to achieve a better result.  This open-ended process was repeated several times for ‘idealised’ purposes. Within this socio-technical context, the Preferendus served as a design/decision support tool to arrive at the best-fit for common-purpose design. To show the real potential of Odesys’ Preferendus, the IMAP results are compared both with single-objective design outcomes and with the min-max goal attainment outcomes.

\subsection{Design application 1: a rail level-crossing service life design}
\textit{Technical context}: Railways and roads often cross each other at level-crossings. Because heavy vehicles must also be able to cross, the railway crossing is often cast in a concrete foundation. The mechanical properties of this concrete foundation are very different from the foundation of the other parts of the railway track. As a result, transitional radiation occurs during the passage of a train, potentially resulting in faster degradation of the local rail system or a negative passenger experience due to vibrational hindrance \citeshortp{wolfert1998,metrikine1998transition}. Therefore, a transition zone is created by varying the number of sleepers and the distance between them to contribute to a smoother transition, which should have a positive effect on both operational performance and passenger comfort. 

\textit{Social context}: In this application, a Multi-Objective Design Optimisation (MODO) of the transition zone is demonstrated, based on several conflicting interests of multiple stakeholders: i.e. (1) capital investment and (2) operational maintenance expenditures, and (3) travel comfort objective functions. It is assumed that these three objectives are linked to three different stakeholders. Take for instance the Dutch ProRail organisation, where there is both a project delivery and a service operations department. They are linked to the capital and operational expenditure objectives, respectively. The Dutch train passenger is represented as the stakeholder linked to the travel comfort objective.

Now, the integrative design problem is firstly described by working through the conceptual threefold framework, see \autoref{fig:threefold}, resulting in design performance, objective and preference functions.

\subsubsection{Design performance functions}
In reality, this design depends on a multitude of design variables, but for now, it will be limited to just two of them:
\begin{enumerate}
    \item $F_1 = x_1 \text{\ } (> 0)$: the distance between the sleepers. Sleepers are the concrete (or sometimes wooden) beams that support the rails, as part of the ballast bed.
    \item $F_2 = x_2 \text{\ } (\ge 1)$: the number of sleepers in the transition zone. The transition zone consists of a different type of sleeper than the rest of the track.
\end{enumerate}

\noindent
Note that (1) in order to be consistent with the general mathematical statement from section 1, the design performance functions $F_1$ and $F_2$ are added here, equal to $x_1$ and $x_2$ respectively, and (2) from the practical application context, the design variables are bounded by $0.3\text{\ }m\le x_1\le0.7\text{\ }m$ and $4\le x_2\le15$, which defines the design space (i.e. the solution space defined by the design variables).

The key design performance functions describing the dynamic behaviour of the track at the level-crossing transition zone are the force $F_3=F(x_1,x_2)$ and the acceleration $F_4=a(x_1,x_2)$. These are usually the result of extensive numerical finite element and/or analytical calculations. For this design application, the physical/mechanical relationships between the design variables are simplified by using interpolation of discrete numerical calculations derived from a finite element based structural dynamic model \citeshortp{shang1}. These interpolated results are the input to the design performance functions.

\subsubsection{Objective functions}
As mentioned before, three objective functions are investigated in this design application: maintenance costs, travel comfort and investment costs. Given these three objectives, the optimal design for the level-crossing zone is determined.\\

\noindent
\textit{Objective maintenance costs (OPEX)}\\
The design of the transition zone is mainly driven by the associated maintenance costs. Large forces and accelerations will have a negative effect on the degradation of the track and foundations, resulting in increased maintenance costs. Hence, this objective can be written as a function of the force and acceleration. For that purpose, the force and acceleration are normalised and combined via the root sum of the square. The final maintenance costs per year objective reads as:

\begin{equation}
    \label{rail:maintenance_costs}
    O_{M}= \sqrt{F_N^2+a_N^2} \cdot 15\,000
\end{equation}

\noindent
where

\begin{equation}
    \label{rail:normalized_force}
    F_N = \frac{F-F_{min}}{F_{max}-F_{min}}
\end{equation}

\begin{equation}
    \label{rail:normalized_acc}
    a_N = \frac{a-a_{min}}{a_{max}-a_{min}}
\end{equation}

\noindent
and where $O_M$ expresses the OPEX per year in EUR.  Note that at the level of design performance functions (i.e. capability-object level), it holds that $F_3 =F$ and $F_4 =a$ respectively.\\

\noindent
\textit{Objective travel comfort}\\
Passenger comfort is an important consideration in railway design. When the dynamic behaviour (due to transition accelerations) during a passage of a level-crossing is substantial, it may lead to a negative travel experience or, in the worst case, to minor mishaps in the train (falling while walking, spilling drinks, etc.). To integrate this into the design problem, an objective is added that describes travel comfort as a function of the normalised acceleration:

\begin{equation}
    \label{rail:travel_comfort}
    O_{C}= 1 - a_N
\end{equation}

\noindent
with $a_N$ as given in Equation (\ref{rail:normalized_acc}).\\

\noindent
\textit{Objective investment costs (CAPEX)}\\
Finally, the investment costs must be considered. The installation of more sleepers will result in higher investment costs. However, more sleepers spread out over a greater distance will also mean that the investment costs for other parts of the rail will be reduced. Therefore, the investment costs objective can be represented as follows:

\begin{equation}
    \label{rail:investment_costs}
    O_{I}= 1000x_2 - 350 x_1 x_2
\end{equation}

\noindent
where $O_I$ expresses the CAPEX in EUR.

\subsubsection{Preference functions}
The preference functions for this design application are constructed based on the input from relevant stakeholders \citeshortp{shang1,shang2021systems}. The three resulting functions, which describe the relations between different values for $P_{1..3,1..3}$ and $O_{1..3}$, are shown as blue curves in \autoref{fig:preference_curves_rail}. Note that the preference function elicitation was performed using the fundamentals of PFM research by \citet{arkesteijn2017improving}.\\ 

\begin{figure}
\centering
\includegraphics[width=\linewidth]{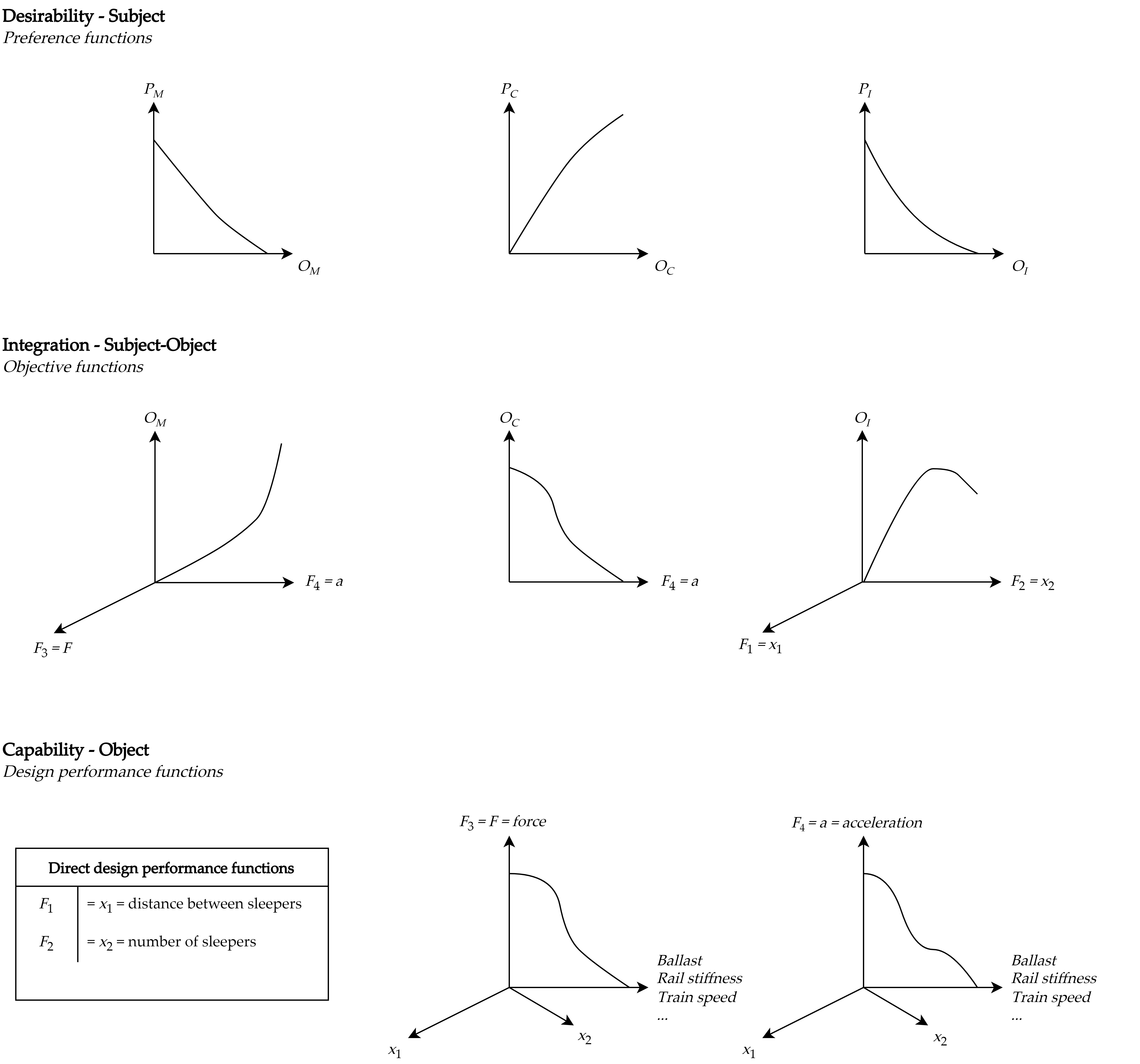}
\caption{Conceptual threefold diagram, describing the systems design integration for the rail level-crossing design application. Note: the aim of this figure is to illustrate the relationship between the different functions and some curves may not represent the actual function.}
\label{fig:threefold_rail}
\end{figure}

\noindent
The systems design integration problem statement is now conceptualised with the threefold diagram in \autoref{fig:threefold_rail}.

\subsubsection{Design optimisation results \& conspection}
To generate the design points (i.e. design configuration results) for the different multi-objective optimisation methods (MODO min-max and IMAP), the weights for each objective must first be determined. Since traditional (contractor) design offices often give a dominant weight to investment costs alone and less to the quality of service (QoS) oriented interests of maintenance and travel performance, here it is deliberately done ‘the other way round’, resulting in $w_{1,M}=0.4$ for maintenance, $w_{2,C}=0.4$ for travel comfort and $w_{3,I}=0.2$ for investments. For evaluation purposes, the design points for the different ($1...3$) Single-Objective Design Optimisations (SODO) are also determined for maintenance-, investment costs and travel comfort respectively.

\begin{figure}
\centering
\includegraphics[width=0.75\linewidth]{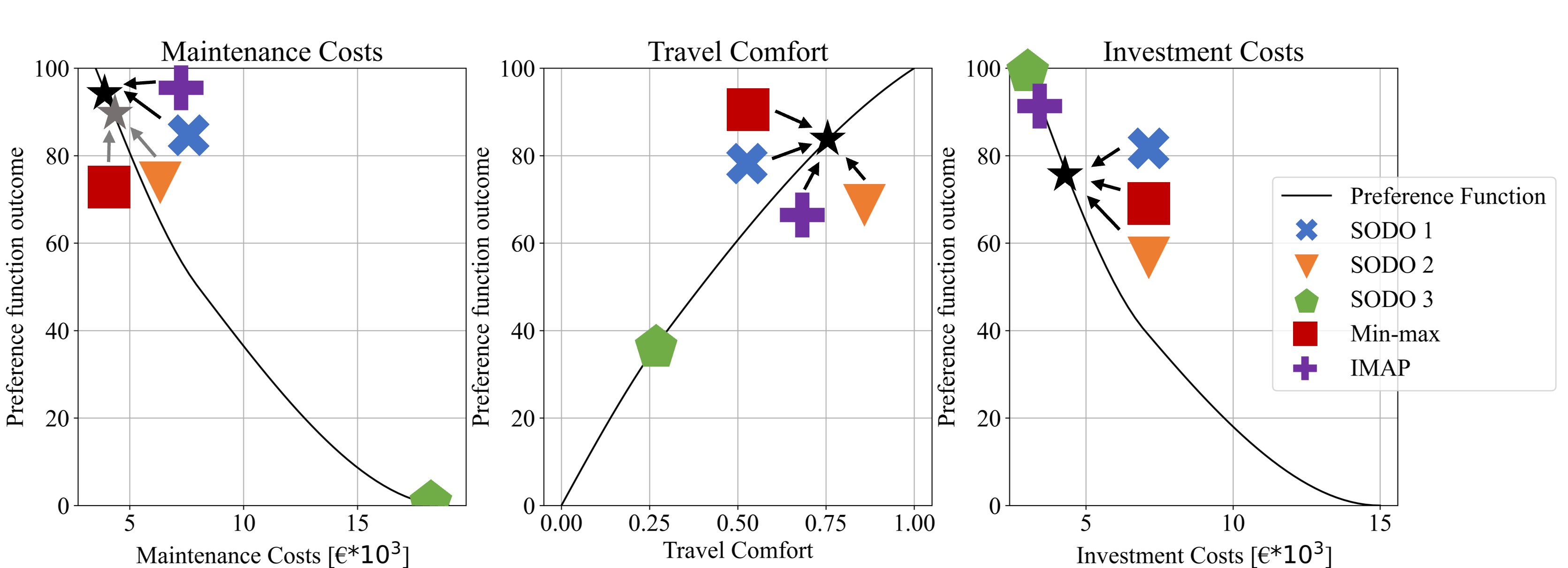}
\caption{The three stakeholder preference functions ($P_{1..3,1..3}$) for different objectives ($O_{1..3}$) for the level-crossing design application, including the results of the different optimisations. The numerical results can be found in \protect\autoref{tab:support_pref_curves_rail}.}
\label{fig:preference_curves_rail}
\end{figure}

The outcomes of the different design points/configurations per optimisation method are first plotted in the preference functions showing the different objective values ($O_{1..3}$) and their corresponding individual preference values ($P_{1..3,1..3}$), see \autoref{fig:preference_curves_rail} and \autoref{tab:support_pref_curves_rail}. Note that, the results for IMAP were obtained with the new Preferendus tool and the other design optimisation results were obtained using specific standard Python routines (see the data availability statement for the repository containing the design application’s code).

Secondly, the numerical results of the different design points/configurations per optimisation method (SODO and/or MODO) can be read from \autoref{tab:eval_rail}. In this table one can find the aggregated preference score, which was used to determine the overall score/ranking via the PFM-based MCDA tool Tetra (the resulting aggregated preference scores are re-scaled between scores of 0 and 100, where 0 reflects the ‘worst’ scoring configuration/alternative and 100 the ‘best’, see \autoref{app:solver} for further details). Note that at least three alternatives are needed for such an overall evaluation (e.g. one reference configuration and two different MODO configurations).

\begin{figure}
\centering
\includegraphics[width=0.6\linewidth]{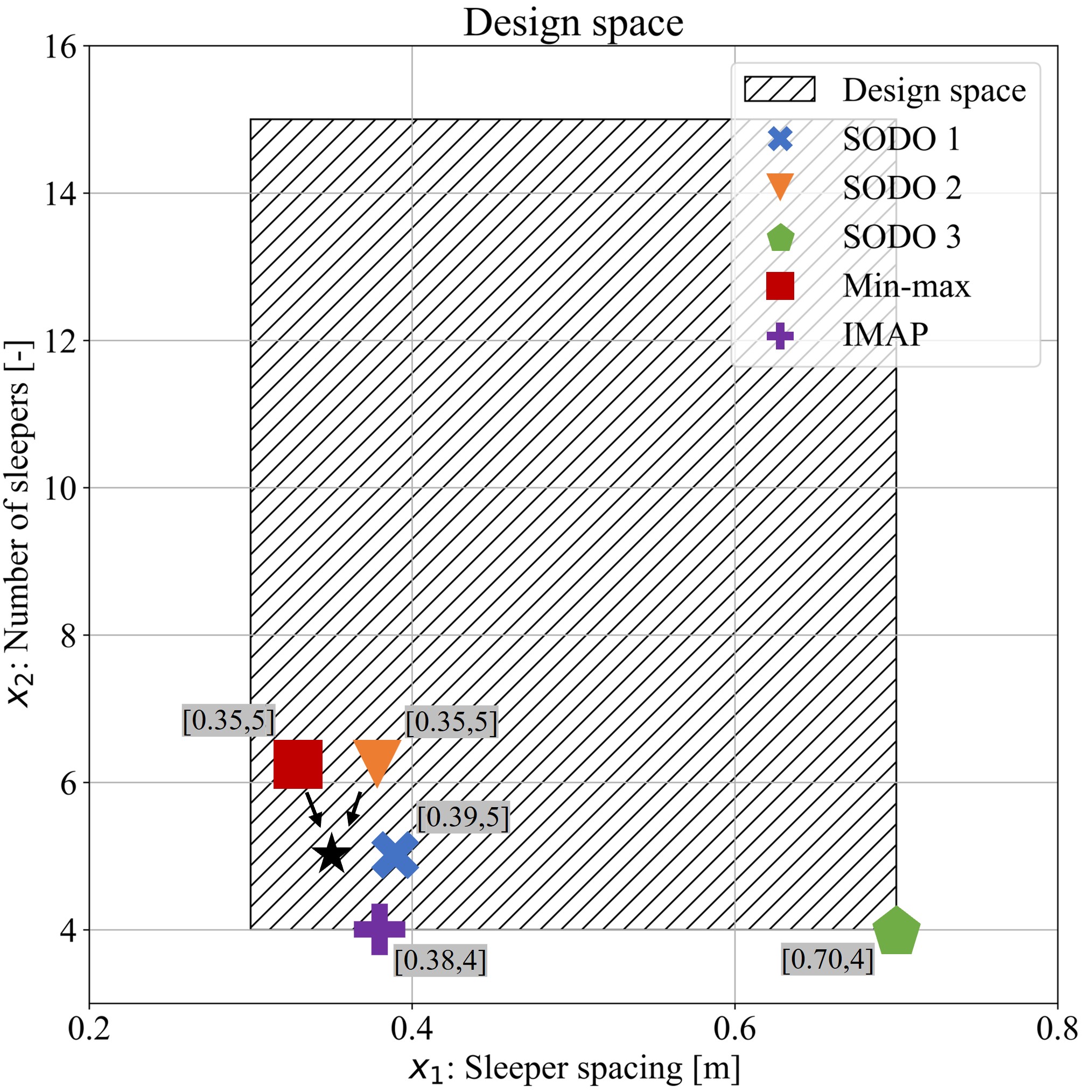}
\caption{The design space of the level-crossing design application and the design configuration/points for the different optimisation methods. The numerical results can be found in \protect\autoref{tab:eval_rail}.}
\label{fig:solution_space_rail}
\end{figure}

Since there are only two design variables in this design application, the two-dimensional design space (sometimes referred to as solution space, see \citet{Dym2004}) containing the design points/configurations per optimisation method can be plotted, see \autoref{fig:solution_space_rail}.\\

\begin{table}[]
\tbl{Results of the objective functions ($O_{1..3}$) and the corresponding preference functions ($P_{1..3,1..3}$) of the level-crossing design application.}
{\begin{tabular}{l c c c c c c}
\toprule
\multicolumn{1}{l}{Optimisation methods} &
\multicolumn{1}{l}{$O_M$ [\texteuro]} &
\multicolumn{1}{l}{$P_M$}  &
\multicolumn{1}{l}{$O_C$} &
\multicolumn{1}{l}{$P_C$}  &
\multicolumn{1}{l}{$O_I$ [\texteuro]} &
\multicolumn{1}{l}{$P_I$} \\
\midrule
Single objective $O_M$ (SODO 1)  &   3942  & 94    &   0.75 & 83   &   4319 & 76    \\
Single objective $O_C$ (SODO 2)  &   4297  & 90    &   0.76 & 84   &   4381 & 75    \\
Single objective $O_I$ (SODO 3)  &   18243 & 0     &   0.27 & 36   &   3020 & 100    \\
MODO min-max                    &   4305  & 90    &   0.76 & 84   &   4382 & 75    \\
MODO IMAP                       &   3974  & 94    &   0.75 & 83   &   3466 & 91    \\
\bottomrule
\end{tabular}}
\label{tab:support_pref_curves_rail}
\end{table}

\begin{table}[]
\tbl{Evaluation of different design configurations per optimisation method and their relative ranking (based on aggregated preference scores) for the level-crossing design application.}
{\begin{tabular}{l l c c}
\toprule
\multicolumn{1}{l}{Optimisation methods} &
\multicolumn{1}{l}{$x_1$ [$m$]} &
\multicolumn{1}{l}{$x_2$}  &
\multicolumn{1}{l}{Aggregated preference score} \\
\midrule
Single objective $O_M$ (SODO 1)  &   0.39    &   5   &   84    \\
Single objective $O_C$ (SODO 2)  &   0.35    &   5   &   81    \\
Single objective $O_I$ (SODO 3)  &   0.70    &   4   &   0    \\
MODO min-max                    &   0.35    &   5   &   81   \\
MODO IMAP                       &   0.38    &   4   &   100   \\
\bottomrule
\end{tabular}}
\label{tab:eval_rail}
\end{table}

\newpage\noindent
The following three conclusions can be drawn from these figures and tables:

(1) The IMAP configuration is either equal to or closest to the best result on all single objectives (the SODO configurations). Only for the single-objective investment costs, IMAP is second best, since it also aims to optimise the other two objectives $O_M$ and $O_C$. For these objectives, a low sleeper spacing ($x_1$) is expected, while the number of sleepers ($x_2$) has a relatively small influence on the outcome of these objectives. For this design application in particular, and given the different objectives and associated stakeholder preferences, a low sleeper spacing ($x_1$) is expected to have a significant impact on objectives $O_M$ and $O_C$, while the number of sleepers ($x_2$) will have a smaller impact.

However, for objective $O_I$, the influence of $x_2$ will be significant, because for lower $x_2$ the investment costs decrease. Furthermore, the influence of $x_1$ on $O_I$ is opposite to its influence on the other two objectives. Therefore, the design configuration that is optimised for investment costs only is not representative. A MODO optimisation is expected to find the ideal balance for the sleeper spacing ($x_1$), with the number of sleepers on the lower bound (i.e. $x_2=4$). The result of the IMAP optimisation does indeed reflect this best-fit for common-purpose balance. As a result, IMAP may be characterised as a pure synthesis multi-objective design method.

(2) The IMAP configuration achieves better or equal individual preference function values ($P_{1..3,1..3}$) and, more importantly, much better overall scores than the MODO min-max method result. This is because, according to the min-max principle, this method will not be able to outperform the one objective score that shows the maximum attainable minimum distance to 100 (i.e. the minimum dissatisfaction). Thus, the min-max method inherently produces a sub-optimal compromise design configuration which, depending on the specific input parameters, can at best perform as well as the synthesis IMAP method. This limits the applicability of the min-max method as a real multi-objective design optimisation method.

(3) From the design space figure it is seen that, perhaps counter-intuitively, both the SODO 1 and 2 and the MODO min-max results fall within the design space ($x_1;x_2$ equals 0.35/0.39 and 5 respectively) and that the MODO IMAP and SODO 3 results lie on the edge and in a corner point of the design space respectively. This is because the set of design points that fall within the design space are the result of optimising the ’technical’ design performance only. In other words, this means that these optimal solutions move to an optimum only within the feasibility space (i.e. a solution space defined by the physical engineering variables only, and which is a subset of the design space) and lie on the classical Pareto front. Note that in this case a possible Pareto front, which defines an edge of the feasibility space as a function of $F$ and $a$, results only from the minimisation of $O_M$ and $O_C$. Despite the fact that SODO 3 actually does find the edges of the design space (corner point), it still scores low overall because it is by far the lowest on the other two objectives (1 and 2). MODO IMAP gives the overall best design point on the edge of the design space ($x_1$ and $x_2$ equal 0.38 and 4 respectively), and can therefore be considered the pure best-fit for common-purpose design point.\\

\noindent
Note that when the emphasis in the design application is on optimising the integrated socio-technical problem, the overall best configuration will be found within and/or on the edge of the design space. When optimising solely on cost or technics, one can either end up at the classical Pareto front or in a corner point of the design space (see also the next design application).

\subsection{Design application 2: a floating wind farm installation design}
\textit{Technical context}: A promising solution for wind energy production in deep waters could be the use of floating wind turbines (FWT). Rather than being placed on a fixed monopile, these turbines are placed on a platform moored to the seabed by anchors. The floating wind farm considered in this design application consists of 36 FWTs and 108 suction anchors (i.e. 3 anchors per FWT).

\textit{Social context}: This application illustrates a MODO approach for the installation of multiple FWTs, taking into account several conflicting interests of multiple stakeholders: i.e. (1) project duration, (2) installation costs, (3) fleet utilisation, and (4) CO2 emissions. Given these four overall interests, an energy service provider (stakeholder one, e.g. Shell) requires a marine contractor (stakeholder two, e.g. Boskalis) to determine the optimal installation design plan. While cost remains a significant factor in the offshore industry, the energy service provider’s primary concern lies in minimising delivery time to expedite resource income generation. Secondly, the energy service provider will have an interest in reducing the CO2 emissions of the project, as this will benefit its carbon footprint and the societal acceptance of the project. The marine contractor’s primary focus will be on reducing the costs, as this will make it more competitive. Secondly, the fleet management department may express a preference for optimising fleet utilisation to maximise operational efficiency.

Now, the integrative problem is firstly described by working through the conceptual threefold framework, see \autoref{fig:threefold}, resulting in design performance, objective, and preference functions.

\subsubsection{Design performance functions}
Several types of vessels are available for the installation of the FWTs and their suction anchors. The amounts of vessels used in the project are the initial three design variables:
\begin{enumerate}
    \item $F_1=x_1 \text{\ } (0\le x_1 \le3)$: small offshore construction vessels (OCV), capable of carrying up to 8 anchors.
    \item $F_2=x_2 \text{\ } (0\le x_2 \le2)$: large offshore construction vessels, capable of carrying up to 12 anchors.
    \item $F_3=x_3 \text{\ } (0\le x_3 \le2)$: self-propelled crane barges, capable of carrying up to 16 anchors.
\end{enumerate}

\noindent
Note that the lower bound of these three design variables is equal to zero. Therefore, a design performance constraint is required to ensure that the sum of all vessels on the project is greater than one (reflecting that at least one vessel is required):

\begin{equation}
    \label{fw:g1}
    g_1 = -(F_1 + F_2 + F_3) + 1 \le 0
\end{equation}

\noindent
This design application also considers the design of the anchors themselves. To do this, design performance functions are defined that describe: (1) the resistance of the anchor to the forces acting on it, and (2) the amplitude of the forces acting on the anchor.

The resistance of the anchors considered in this design application can be estimated using analytical design calculations according to \citet{houlsby2005design,randolph2017offshore,arany2018simplified}. These calculations usually depend on several design variables, only two of which are considered here:

\begin{enumerate}
    \item $F_4=x_4 \text{\ } (> 0)$: Diameter of the suction anchor in meters.
    \item $F_5=x_5 \text{\ } (> 0)$: Penetration length of the suction anchor in meters.
\end{enumerate}

\noindent
For practical reasons, these variables are bounded by $1.5\text{\ }m\le x_4\le4\text{\ }m$ and $2\text{\ }m\le x_5\le8\text{\ }m$. The other design variables are uncontrollable variables $\mathbf{y}$ in this design application, where $\mathbf{y}=$ [working point $F_a$, mooring configuration, anchor type, soil conditions, mooring line properties]. Consequently, the anchor resistance can be mathematically formulated as $F_6=R_a\left(x_4,x_5,\mathbf{y} \right)$. The soil is assumed to be clay with an undrained shear strength of $s_u=60 \text{\ } kPa$ and a submerged weight of $\gamma '=9 \text{\ } kN/m^3$. The coefficient of friction between the anchor shaft and the soil is $\alpha=0.64$. The mooring line consists entirely of a chain with a nominal diameter of $240\text{\ }mm$. This chain is attached to the anchor at a depth of 0.5 times the penetration length. Furthermore, the coefficient of friction between the seabed and the chain is taken as $\mu=0.25$ and the active bearing area coefficient \emph{AWB} $=2.5$.

While anchor resistance can be determined by analytical calculations, the forces acting on the anchor cannot be determined in the same manner. This is due to their dependence on various variables such as platform type, mooring line characteristics, pre-tension, and anchor radius. To obtain accurate normative forces, numerous numerical time-domain calculations must be performed, as outlined in \citet{dnv_0437}. These calculations are beyond the scope of this paper. Instead, the relevant design variables are considered as uncontrollable physical variables $\mathbf{y}$, resulting in the following (assumed) force on the anchors: $F_7=F_a\left(\mathbf{y}\right)=3.8 \text{\ } {M\!N}$, where $\mathbf{y}=$ [platform type, mooring line characteristics, pre-tension, mooring line length, anchor radius].

The two design performance functions $F_6$ and $F_7$ are related through a design performance constraint. This constraint describes (part of) the feasibility space of the ‘technical’ design by defining the boundary where the resistance of the anchor is greater than or equal to the force on the anchor:

\begin{equation}
    \label{fw:g2}
    g_2 = F_7(\mathbf{y}) - F_6(x_4,x_5,\mathbf{y}) = F_a - R_a \le 0
\end{equation}

\subsubsection{Objective functions}
As mentioned before, four objectives are investigated in this design application: project duration, installation costs, fleet utilisation, and CO\textsubscript{2} emissions. Given these four objectives, the optimal design plan for installing the FWTs is determined.\\

\noindent
\textit{Objective project duration}\\
The project duration depends on the number of vessels involved in the project, their deck capacity and the speed at which they can install anchors, which is assumed to be one anchor/day/vessel. In addition, after all the anchors on board have been installed, the vessels will have to load new anchors. This process takes 1.5 days for the small OCV, 2 days for the large OCV, and 2.5 days for the barge.

To obtain the overall project duration, a discrete event simulation (DES) was incorporated into the model, which depends on the type and number of vessels (i.e. $x_1..x_3$). See the data availability statement for the code of the DES. In conclusion, the objective function for the project duration can be expressed as follows:

\begin{equation}
    \label{fw:proj_duration}
    O_{PD}=f(x_1,x_2,x_3)
\end{equation}

\noindent
where $f$ is the DES and $O_{PD}$ is expressed in days.\\

\noindent
\textit{Objective installation costs}\\
The project’s installation costs objective depends on two components: (1) the day rates of the vessels, and (2) the cost of the anchors. The following theoretical day rates $R$ are assumed:
\begin{enumerate}
    \item Small OCV ($x_1$): $R_1 = \text{\texteuro} 47,000/day$
    \item Large OCV ($x_2$): $R_2 = \text{\texteuro} 55,000/day$
    \item Barge ($x_3$): $R_3 = \text{\texteuro} 35,000/day$
\end{enumerate}

\noindent
The cost per anchor can be divided into a fixed part (\texteuro 40,000/anchor) and a variable part, where the variable part depends on the material costs (\texteuro815/t). This results in the following objective cost function:

\begin{equation}
    \label{fw:costs}
    O_{C}=(815M_a+40,000)n_a+\sum_{i=1}^3 x_i t_i R_i
\end{equation}

\noindent
where $O_C$ is expressed in EUR, $n_a$ is the number of anchors (i.e. $n_a=108$), $t_i$ the time a vessel is needed (result from the DES), and $M_a$ the mass of the anchors, which is defined as:

\begin{equation}
    \label{fw:mass_anchor}
    M_a=\left(\pi x_5x_4t+\frac{\pi}{4}x_4^2t \right)W_{steel}
\end{equation}

\noindent
with $W_{steel}$ is the weight of steel, assumed as $78.5 \text{\ } t$ (‘tonnes’).\\

\noindent
\textit{Objective fleet utilisation}\\
For a marine contractor, optimal fleet utilisation is a key driver. Consequently, this objective focuses on evaluating the probability of a vessel being better utilised in another project (e.g. specialised vessels are preferred to multi-purpose vessels). For this purpose, the following values are assumed:

\begin{enumerate}
    \item Small OCV ($x_1$): $p_1 = 0.7$
    \item Large OCV ($x_2$): $p_2 = 0.8$
    \item Barge ($x_3$): $p_3 = 0.5$
\end{enumerate}

\noindent
The fleet utilisation objective is then defined as:

\begin{equation}
    \label{eq:fleet_util}
    O_{F}= \prod\limits_{i=1}^3 p_i^{x_i}
\end{equation}

\noindent
where $O_F$ is expressed as the combined chance with a value between [0, 1].\\

\noindent
\textit{Objective CO\textsubscript{2} emissions}\\
Sustainability is becoming an increasingly important aspect within offshore (wind) project development. Most of the emissions will be generated by the vessels, for which the following theoretical average emission rates are assumed:

\begin{enumerate}
    \item Small OCV ($x_1$): $E_1 = 30 \text{\ }t/day$
    \item Large OCV ($x_2$): $E_2 = 40 \text{\ } t/day$
    \item Barge ($x_3$): $E_3 = 35 \text{\ } t/day$
\end{enumerate}

\noindent
As other sources of emissions are neglected, the emission objective is defined as:

\begin{equation}
    \label{fw:sustainability}
    O_{S}=\sum_{i=1}^3 x_iE_it_i
\end{equation}

\noindent
where $O_S$ is expressed in tonnes and with $t_i$ the time a vessel is needed (result from the DES).\\

\noindent
Note that the Odesys mathematical statement allows for the direct integration of design performance and objective functions. However, in certain cases, design performance functions will not only directly link to the objective functions but can also connect through (in)equality design performance constraints. This indirect linking is common in design problems where, for example, force constraints play an important role. In such cases, these constraints define the feasibility space, and together with directly linked design performance functions, they span the design (i.e. solution) space.

\subsubsection{Preference functions}
The preference functions for this design application were developed with floating wind project experts within Boskalis, based on the input from an energy service provider. The four resulting functions, which describe the relations between different values for $P_{1..2,1..4}$ and $O_{1..4}$), are shown as blue curves in \autoref{fig:preference_curves_fw}. Note that the
preference function elicitation was again (like in the previous design application) performed using the fundamentals of PFM research by \citet{arkesteijn2017improving}.\\ 

\begin{figure}
\centering
\includegraphics[width=\linewidth]{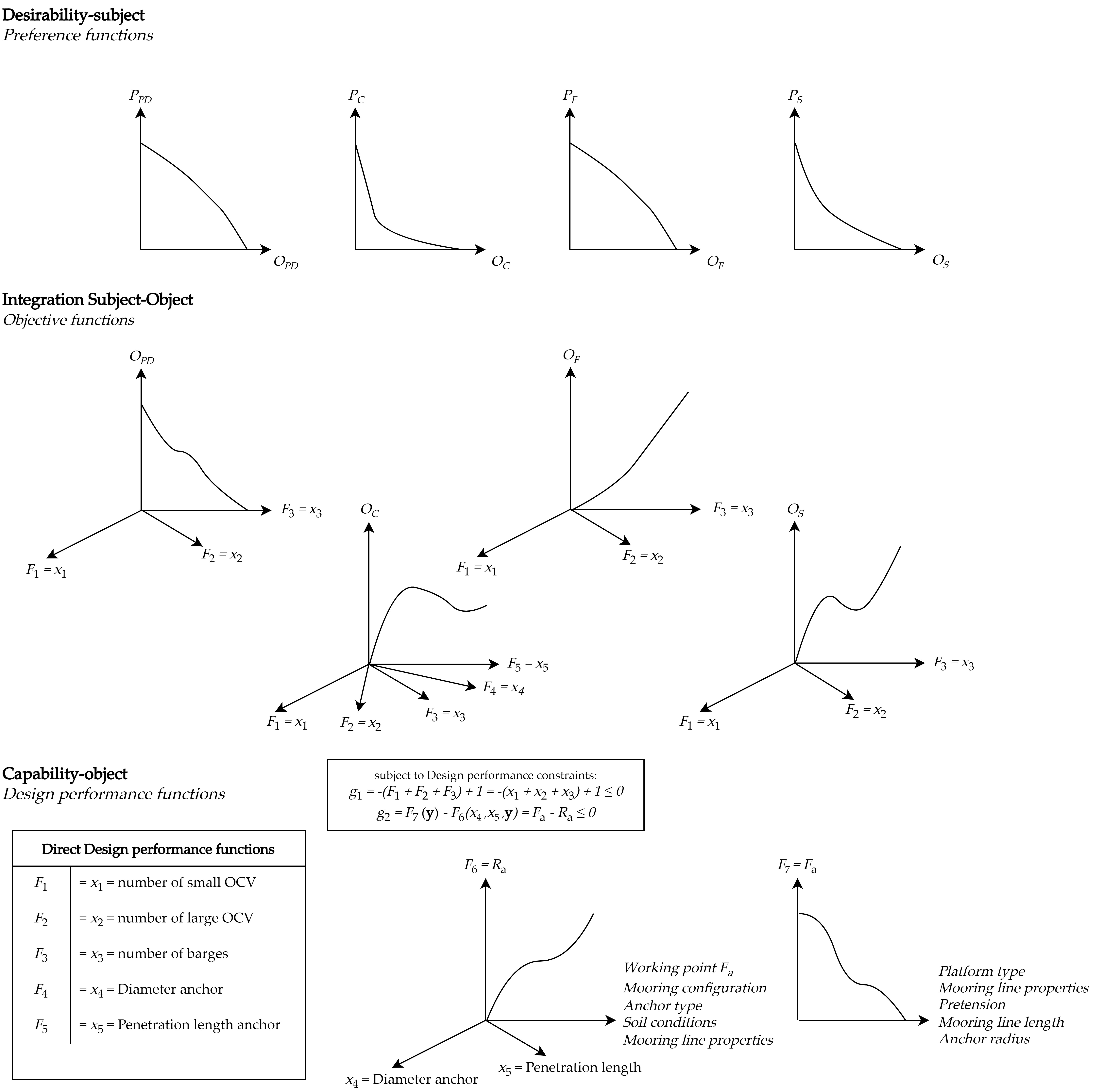}
\caption{Conceptual threefold diagram, describing the systems design integration for the floating wind turbine design application. Note: the aim of this figure is to illustrate the relationship between the different functions and some curves may not represent the actual function.}
\label{fig:threefold_fw}
\end{figure}

\noindent
The systems design integration problem statement is now conceptualised with the
threefold diagram shown in \autoref{fig:threefold_fw}.

\subsubsection{Design optimisation results \& conspection}
To generate the design points (i.e. design configuration results) for the different multi-objective optimisation methods (MODO min-max and IMAP), the weights for each objective must first be determined. Traditionally, installation costs have been the main driver for offshore projects and/or tender bids. However, with the introduction of the Odesys design optimisation methodology, it is now possible to optimise the design considering other relevant objectives that reflect the shared value of the installation plan for both the energy service provider and the contractor. The following weight distributions were chosen to model this joint plan: $w_{1,PD}=0.30$ for project duration, $w_{1,S}=0.20$ for sustainability (emissions), $w_{2,C}=0.35$ for the installation costs, and $w_{2,F}=0.15$ for fleet utilisation.

For evaluation purposes, both the single-objective optimisation of $O_C$ (SODO costs) and the MODO min-max optimisation design points are also determined. Note that the other SODOs (single-objective optimisations on $O_{PD}$, $O_F$, and $O_S$) cannot be included in the integral evaluation as they are not dependent on $x_4$ and $x_5$ (but only on $x_1..x_3$). 

\begin{figure}
\centering
\includegraphics[width=\linewidth]{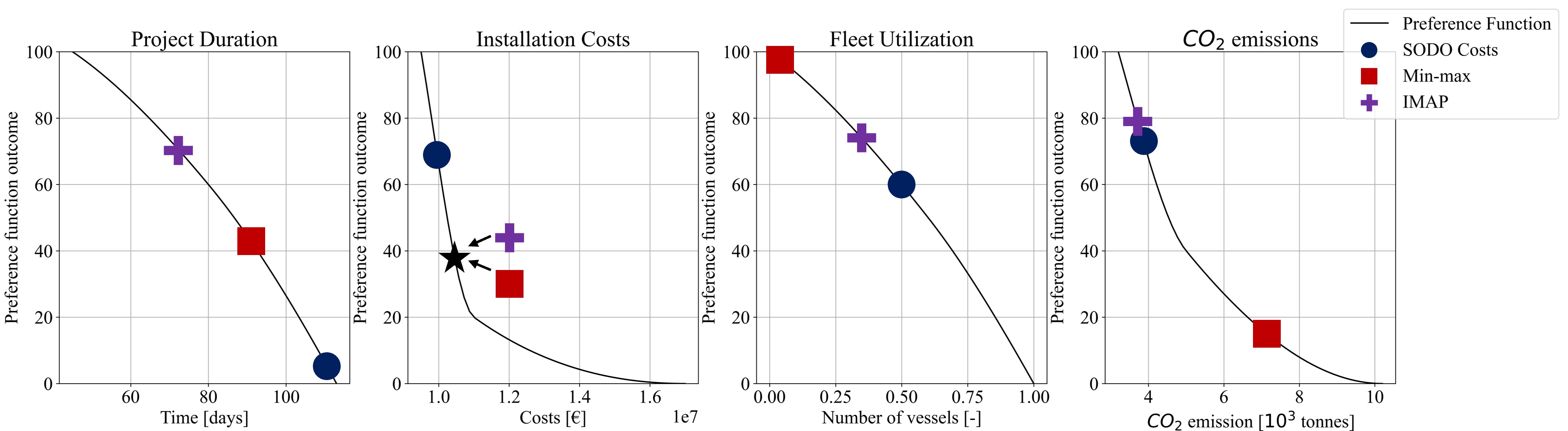}
\caption{The four stakeholder preference functions ($P_{1..4,1..4}$) for different objectives ($O_{1..4}$) for the floating wind design application, including the results of the different optimisations. The numerical results can be found in \protect\autoref{tab:support_pref_curves_fw}.}
\label{fig:preference_curves_fw}
\end{figure}

The outcomes of the different design points/configurations per optimisation method are first plotted in the different preference functions showing the different objective function values ($O_{1..4}$) and their corresponding individual preference function values ($P_{1..2,1..4}$), see \autoref{fig:preference_curves_fw} and \autoref{tab:support_pref_curves_fw}. Secondly, the numerical results of the different design points/configurations per optimisation method can be read from \autoref{tab:eval_fw}. In this table, one can also find the aggregated preference score, which was used to determine the overall score/ranking via the PFM-based MCDA tool Tetra (the resulting aggregated preference scores are re-scaled between scores of 0 and 100, where 0 reflects the ‘worst’ scoring configuration/alternative and 100 the ‘best’, see \autoref{app:solver} for further details).\\

\begin{table}[]
\tbl{Results of the objective functions ($O_{1..4}$) and the corresponding preference functions ($P_{1..4,1..4}$) of the floating wind design application.}
{\begin{tabular}{l c c c c c c c c}
\toprule
\multicolumn{1}{l}{Optimisation methods} &
\multicolumn{1}{l}{$O_{PD}$ [$days$]} &
\multicolumn{1}{l}{$P_{PD}$}  &
\multicolumn{1}{l}{$O_C$} [\texteuro] &
\multicolumn{1}{l}{$P_C$}  &
\multicolumn{1}{l}{$O_F$} &
\multicolumn{1}{l}{$P_F$} &
\multicolumn{1}{l}{$O_S$ [$t$]} &
\multicolumn{1}{l}{$P_S$} \\
\midrule
Single objective $O_C$ (SODO costs)  &   110.5 & 5     &   9.96E6  & 69   &   0.50 & 60   & 3868 & 73  \\
MODO min-max                         &   91    & 43    &   10.45E6 & 38   &   0.04 & 97   & 7135 & 15  \\
MODO IMAP                            &   72.5  & 70    &   10.47E6 & 37   &   0.35 & 74   & 3722 & 79  \\
\bottomrule
\end{tabular}}
\label{tab:support_pref_curves_fw}
\end{table}

\begin{table}[]
\tbl{Evaluation of different design configurations per optimisation method and their relative ranking (based on aggregated preference scores) for the floating wind design application.}
{\begin{tabular}{l l c c c c c}
\toprule
\multicolumn{1}{l}{Optimisation methods} &
\multicolumn{1}{l}{$x_1$} &
\multicolumn{1}{l}{$x_2$} &
\multicolumn{1}{l}{$x_3$} &
\multicolumn{1}{l}{$x_4$ [$m$]} &
\multicolumn{1}{l}{$x_5$ [$m$]} &
\multicolumn{1}{l}{Aggregated preference score} \\
\midrule
Single objective $O_C$ (SODO costs)  &  0  &  0  &  1  & 2.2  & 8.0  &  69  \\
MODO min-max                    &  1  &  0  &  2  & 2.2  & 8.0  &  0  \\
MODO IMAP                       &  1  &  0  &  1  & 2.2  & 8.0  &  100   \\
\bottomrule
\end{tabular}}
\label{tab:eval_fw}
\end{table}

\noindent
The following three conclusions can be drawn from these figures and tables:

(1) Comparing the IMAP configuration with the SODO design point on installation costs, IMAP outperforms the SODO on three of the four objectives. This difference is most evident when the result of the project duration objective is compared with the result of the installation cost objective. These objectives are opposite by the impact of the number of vessels ($x_{1..3}$) on them.  More vessels lead to faster project completion but higher costs. Therefore, a design configuration that scores well on cost will not score well on project duration, as can be seen for the SODO on installation costs. This result illustrates that considering cost alone (single stakeholder and single objective approach) is not an accurate reflection of the real planning challenge. In contrast, IMAP demonstrates a balanced approach by considering multiple objectives, including both the technical design and economics.

(2) The overall score of the IMAP configuration is substantially higher than that of the min-max method. As the min-max method tries to minimise the distance to a score of 100 for all different preference scores $P_{1..2,1..4}$, it can result in very low preference scores for conflicting objectives. In this design application, this is the case for the project duration ($O_{PD}$) and installation costs ($O_C$) objectives. As a result, the min-max solution scores low for these two objectives. This is in contrast with the IMAP design solution, which can find higher preference scores $P_{1..2,1..4}$ for these two objectives. The presence of these conflicting interests thus limits the applicability of the min-max method, as also shown in the first design application. Note that it can still perform well for a ‘single’ interest, as shown by the positive reflection of the fleet utilisation objective with the use of more barges.

(3) \autoref{tab:eval_fw} shows that all three solutions have the same result for design variables $x_4$ and $x_5$. This indicates that this particular combination of $x_4$ and $x_5$ yields the lowest anchor cost without violating the design performance constraint $g_2$. In other words, for all three methods, there would be no difference in the optimisation if limited to a purely technical optimisation within the feasibility space. However, the added value of IMAP is evident from the results for design variables $x_1$, $x_2$, and $x_3$, where IMAP can arrive at an overall better design solution than the other two methods by including both technical and vessel-related installation planning concerns. Note that also the best outcome within the feasibility space for $x_4$ and $x_5$ will change if objectives in the managerial (subject desirability) domain favour technical over dimensioning of the suction anchors. In such cases, the solution may be selected away from the edge of the feasibility space (i.e. the Pareto front) as it offers greater benefits to the overall planning and design performance.

\section{Discussion \& further developments}
Although both design applications are simplified for methodological illustration purposes, they already demonstrate the added value of the Preferendus/IMAP in the field of multi-objective design optimisation. In addition, the Preferendus/IMAP demonstrates its practical value through its application and validation in the following real-life projects: (1) the primary design and construction/production management processes of the marine contractor Boskalis \citep{vanheukelumIALCCE}; (2) the EU NRG-Storage research project \citep{zhilyaev2022best}; (3) several PhD/MSc thesis project applications \citep{shang2021systems,shang1,van2022preference}. In all projects, stakeholders with decision-making power (on both the developer and contractor side) are predominantly positive about the unexpected design solutions they could not have achieved without the use of this computer-aided design ‘engine’, the Preferendus, as part of the Odesys methodology.

For both design applications (rail level-crossing and floating wind turbine), a major extension of the models is currently underway to better fit the design/decision problem in practice, so that both more realistic design performance and also better preference functions will be included. For the floating wind application, this means that OpenFAST, an open-source wind turbine simulation tool, is linked via a surrogate model and integrated at the level of design performance functions. For the level-crossing application, the modelling input will be refined at all levels (focus on the preference and objective functions). In addition, for the floating wind application, but also for a dredging application, validation sessions will be carried out to refine the modelling inputs (especially the performance functions) and to evaluate the results, especially of the new IMAP and the existing min-max methods. This is done in the form of a serious game, using the Preferendus as design support ‘engine’, with the aim of increasing the internal acceptance and the link with the iterative group design engineering process.

Based on these current developments, at least three main focuses for the further development of the Preferendus can be formulated:
\begin{enumerate}
    \item The output depends not only on the best possible design performance functions, but also on a good reflection of human objectives and preferences. Especially for the latter, further preference elicitation research is needed to arrive at balanced preference functions with corresponding individual preferences as input.
    \item The result of the optimisation may also be an empty design solution space: i.e. a so-called stalemate situation. In this case, additional decision support functionality will need to be provided to support stakeholders to achieve the best possible negotiation to still arrive at an acceptable design solution space.
    \item The design performance functions are currently deterministic. However, for more realistic applications, probabilistic design modelling techniques will need to be integrated, e.g. for the offshore design application, uncertainty in working hours or operational weather slots. Improvements to the current Discrete Event Simulator (DES) may be required, particularly for repetitive production and installation operations.
\end{enumerate}

In addition, the Preferendus and the IMAP method will be applied in other systems design and management applications, such as dynamic preference and performance-based mitigation control (MitC) of large construction projects and optimal planning of flood defence system reinforcements (see e.g. \citeshortt{kammouh2022dynamic,klerk2021optimal}). For application within the European NRG-Storage project (European Commission NRG-STORAGE project (no. GA 870114)), the current Preferendus model will be evaluated against the min-max optimisation approach in a real-life context. This will be done as an extension of an MSc thesis project in which the added value of the Preferendus within the so-called social cycle of an urban planning project was studied in more detail within a municipality \citep{van2022preference}.

Furthermore, the added value within the so-called concurrent engineering and design developments in the field of ’Early Contractor Involvement’ is also investigated. In particular, the Preferendus will be used to support and evaluate the new so-called two-phase contract for infrastructure projects, in which the activities of the Dutch national infrastructure service provider (RWS) and its contractors are further intertwined, in order to avoid major contract changes that are the result of the classic serial, non-participative design and engineering process.

Finally, the Odesys methodology has already been taught, and further tested and validated in several MSc courses in Systems Engineering Design at the Faculty of Civil Engineering \& Geosciences at Delft University of Technology this year. The purpose is to further explore the added value and potential improvements of the Preferendus as soon as possible. Within these courses, MSc students develop a Preferendus/IMAP-based model of a self-chosen real-life system of interest as part of the so-called Open Design Learning (ODL) response (for more details on the ODL response, see \citeshortt{wolfert2022fit}). Some findings from these courses have already been incorporated into the current Preferendus code, see the data availability statement for further details.

\section{Conclusions}
The aim of the Odesys methodology is to promote the adoption of engineering artefacts in our future society by following an open space/source design and systems integration approach supported by sound mathematical open glass box optimisation models. This as a means to achieve well informed decision-making leading to the best-fit open-ended solutions for socio-eco purpose. This requires systems thinking and a stakeholder-oriented focus to explore different solutions within an open-ended optimisation process, uniting both capability (technological) derived from the engineering asset’s performance and desirability (sociological) derived from each stakeholder’s preferences. This results in an open dialogue and a co-design approach that enables \textit{a priori} best-fit for common-purpose design synthesis dissolutions rather than \textit{a posteriori} design compromise absolutions.

This paper presents a pure \textit{a priori} socio-technical systems integration and design methodology, together with a new Integrative Maximised Aggregated Preference (IMAP) optimisation method. Furthermore, IMAP has been integrated into the Preferendus tool, which combines state-of-the-art principles of PFM with a specifically developed inter-generational GA optimisation solver. Two specific engineering systems design applications have been worked out by first using the threefold diagram to formulate the mathematical problem statement. The resulting outcomes of these applications clearly demonstrate the added value of IMAP/Preferendus.

Firstly, it provides a single best-fit for common-purpose design point, unlike a Pareto front where a systems designer still has to choose the final design because the front does not define a single optimal design point. This solves an important modelling error, in addition to the fact that classical design optimisation methods leading to these Pareto fronts contain fundamental aggregation errors, namely that design configurations lying on the Pareto front cannot all have the same preference scores.

Secondly, IMAP returns the best design configuration in all design applications compared to a set of single-objective design configurations and a design configuration obtained by the classical multi-objective min-max method. This allows IMAP to be characterised as a real synthesis multi-objective design method that ensures a best-fit for common-purpose point within the design space, rather than a sub-optimal one-sided corner point and/or a best point in the feasibility space only.

Finally, IMAP/Preferendus truly unites design performance functions (supply), via the level of inter-play objective functions, with stakeholder’s preference functions (demand), synthesising for the best-fit for common-purpose solution and outperforming one-sided design approaches that focus only on the technical domain. This means that the IMAP/Preferendus is either equal to other design methodologies in the technical domain, but outperforms methodologies within the management domain (see design application 2: a floating wind turbine installation) or outperforms other design methodologies in both the technical and the management domains (see design application 1: a rail level-crossing service life design).

\section*{Acknowledgements}
Thanks to Ms Shang (PhD student) for providing the necessary data for the railway infrastructure design application. Thanks to Mr Zhilyaev (PhD student within the EU, NRG-Storage project) for the input on and verification of the mathematical statement. Both are currently PhD students in the Engineering Asset Management group at Delft University of Technology. Finally, thanks to Lukas Teuber and Matt Julseth for proofreading (parts of) this article.\\

\noindent
Part of this paper has been published \textit{in verbatim} in the book Open Design Systems \citep{opendesignsystems}, which provides a more comprehensive embedding and context for the Open Design Systems methodology introduced in this paper.

\section*{Disclosure statement}
The authors report that there are no competing interests to declare.

\section*{Data availability statement}
The Preferendus software tool, including the design applications discussed in this paper, is available via GitHub: \url{https://github.com/TUDelft-Odesys/Preferendus}. Note that the Preferendus application uses an API endpoint on the Tetra server to interact directly with the Tetra solver. This API endpoint requires authentication. For assessment purposes, credentials can be requested via the corresponding author. An interactive version of the Tetra solver can be found at \url{http://choicerobot.com}.

\bibliographystyle{apacite}
\bibliography{interactapasample}

\appendix
\section{The inter-generational GA solver}
\label{app:solver}
To find the design configuration that reflects the integrative maximum aggregation of reference, it is necessary to use an optimisation algorithm. Furthermore, this algorithm must also be able to interact with the Tetra, the multi-criteria decision analysis (MCDA) software tool based on Preference Function Modelling (PFM). The algorithm of the non-linear Tetra solver is based on minimising the least-squares difference between the overall preference score and each of the individual scores (on all decision criteria) by computing its closest counterpart (for more information on the Tetra software, see \citet{scientific_metrics}).

For this purpose, a Genetic Algorithm (GA) has been developed that is specifically tailored to work with a PFM solver (like Tetra) and its specific features of normalised scores and relative ranking. First, these features are described.

\subsection{Normalised scores}
Preference scores are expressed as numbers on a defined scale, here ranging from 0 to 100, where 0 is the ‘worst’ scoring design configuration/alternative and 100 is the ‘best’. This means that the best alternative will always have a score of 100 and the worst alternative will always have a score of 0. Since a GA typically checks whether the best score of the current generation ($G_n$) outperforms the previous one ($G_{n-1}$), normalised scores will lead to convergence problems; the GA cannot determine whether an improvement occurs, since the best alternative will always have a score of 100.

Also, in the case of constrained problems, where the alternative with a score of 100 may be infeasible and should be disregarded, convergence problems persist. In this case, it is possible that the best feasible design alternative will have a lower preference score in generation $G_n$ compared to generation $G_{n-1}$. This is because, due to normalisation, the score of one alternative always depends on the performance of all other alternatives. This must be taken into account in the GA solver.

\subsection{Rank reversal}
Rank reversal, the notion that ranks might change when an alternative is added or removed from the population, is common in various MCDA models and is also present in Tetra \citep{wang2009rank,aires2018rank}. This phenomenon is commonly observed when non-competitive (i.e. irrelevant) alternatives are added or removed from the population \citep{aires2018rank}. In short, especially when extreme or ’irrelevant’ (i.e. without any real meaning) alternatives are added/removed, rank reversal can occur, potentially leading to convergence problems in finding the best solution by evaluating whether the generation ($G_n$) outperforms the previous one ($G_{n-1}$). Furthermore, since an initial population is generated (quasi) randomly, it is not unlikely that extreme or irrelevant alternatives will be part of the first generation evaluated by the GA. These alternatives would never be considered in reality, creating a discrepancy between the GA solver and real-life design alternatives that should be mitigated to achieve convergence.

\subsection{Modifications to the GA}
To solve the aforementioned issues arising from normalisation and rank reversal, the following modifications were applied, resulting in a so-called inter-generational GA solver:

(1) An additional step must be added to the evaluation of a generation. After determining the aggregated preference scores for the entire population, the member with the highest rank is added to a list. This list contains the best members of all generations ($G_n,G_{n-1},...,G_0$) and is evaluated separately to obtain an aggregated preference score for all members in this list. If the aggregated preference score associated with generation $G_n$ is lower than that of $G_{n-1}$, no improvements are made. However, if the score of generation $G_n$ equals is 100, the GA has either improved or, if the score of generation $G_{n-1}$ is also 100, has found a temporary optimum.

(2) The initial population can be built from user-defined initialised solutions. These solutions can be chosen arbitrarily or guided by the single objective and/or min-max design optimisation results. this means that the initial population is no longer (quasi) random and reflects true potential design points, reducing the probability of non-convergence from the start. After the first (initial) population has been evaluated, mutation will begin to diversify the population during the creation of the next generation, making it possible to reach another optimal solution even though the initial population is directionally determined.

Note that this implementation of ’arbitrary’ initialised solutions is also very useful for validating the results. Running the same problem with different starting points can confirm that the result is indeed optimal.

(3) After the initial evaluation of the function $U$ (see Equation (\ref{eq:general_MS})), an additional re-evaluation is introduced by feeding the GA with as many potential real-life design points as possible. This is done by re-evaluating the population in such a way that the very worst alternatives that (potentially) reflect irrelevant, non-competitive alternatives are excluded. This means that after the population has been evaluated for a first time, only alternatives with an aggregated preference score higher than a certain lower limit $P^*$ (which can be set by the designer, here fixed at 20) are re-evaluated a second time, thus improving GA convergence.\\

\noindent
These three modifications have been added to a fit-for-purpose inter-generational solver GA, incorporating key elements from standard available GA Python packages, allowing the aggregated results of one generation to be compared with another. See the data availability statement for details on how to access the code of this solver.

Note that the aforementioned modifications are the result of pragmatic engineering judgement using the principle of reflection and after validation of a large number of example problems. As a possible specific step for further research, it may be of interest (also from the perspective of improved solving speed) to investigate whether other optimisation algorithms than a GA might be more suitable for this specific purpose.

\end{document}